\documentclass{llncs}

\usepackage{a4wide}
\usepackage{amssymb,amsmath}
\usepackage{comment}
\usepackage{complexity}

\usepackage{enumitem}
\usepackage{bm}

\usepackage{graphicx}

\usepackage[algosection,algoruled,vlined,nofillcomment]{algorithm2e}
\SetKw{nil}{NIL}

\usepackage[colorlinks=true,citecolor=blue,urlcolor=blue,linkcolor=blue,bookmarksopen=true,unicode=true,pdffitwindow=true]{hyperref}

\newcommand{\set}[1]{\ensuremath{\left\{#1 \right\}}}

\newcommand{\floor}[1]{\ensuremath{\left\lfloor #1 \right\rfloor}}

\newcommand{\scadical}{\textsc{CaDiCaL}}
\newcommand{\smathcheck}{\textsc{MathCheck}}
\newcommand{\ssat}{\textsc{Sat}}
\newcommand{\sslow}{\textsc{Low}}
\newcommand{\ssmid}{\textsc{Mid}}
\newcommand{\ssup}{\textsc{Up}}
\newcommand{\sslm}{\textsc{LM}}

\newcommand{\cp}{\square}

\newcommand{\defi}[1]{\emph{#1}}

\spnewtheorem{observation}[theorem]{Observation}{\bfseries}{\itshape}

\title{
  Matchings of five directions in hypercube extends to Hamilton cycles and paths with prescribed ends
  \thanks{
Supported by the~Czech Science Foundation grant GA 22-15272S.
  }
}
\author{
  Ji\v{r}\'i Fink\inst{1} \and
  Vojt\v{e}ch Hotmar\inst{1}
}

\institute{
  Department of Theoretical Computer Science and Mathematical Logic \\ 
  Faculty of Mathematics and Physics \\ 
  Charles University in Prague \\
  \email{fink@ktiml.mff.cuni.cz}
}

\begin{document}
\maketitle
\begin{abstract}
	The~$n$-dimensional hypercube graph~$Q_n$ has as vertices all subsets of $\set{1, \ldots, n}$, and an~edge between any two sets that differ in a~single element.
	The~Ruskey-Savage conjecture states that every matching of the~$n$-dimensional hypercube~$Q_n$ can be extended into a~Hamilton cycle.
	We prove that matchings of~$Q_n$ containing edges spanning at most~$d = 5$ directions can be extended into a~Hamilton cycle.
	We also characterize when these matchings of most~$d = 5$ directions can be extended into a~Hamilton path between two prescribed vertices.
	Our proofs work for arbitrary $d$ and $n$ where $d \le n$ assuming some extension properties hold in~$Q_d$ which we verified by a~computer for~$d=5$.
\end{abstract}

\section{Introduction}

The~$n$-dimensional hypercube~$Q_n$ is the~graph whose vertices are all subsets of~$[n]:=\set{1,\ldots,n}$ and whose edges connect sets that differ in a~single element.
There is a~large literature on structural properties of this class of graphs which comes from research on the~topological structure and analysis of hypercubic interconnection networks \cite{xu2013topological}.
For example, Gros in 1872~\cite{Gros} proved that~$Q_n$ contains a~Hamilton cycle for every $n \ge 2$.
A~Hamilton path in the~hypercube~$Q_n$ is in the~literature also called a~Binary Gray Code which is a~listing of all binary strings of length $n$ such that every two consecutive strings differ in exactly one positions.
It is named after Frank Gray who in 1953 constructed a~special Hamilton cycle in~$Q_n$, called Reflected Binary Gray Code, and use it in his U.S. Patent 2632058 for signal processing.
Then, Gray codes have found applications in such diverse areas as information retrieval, image processing or data compression, and alternative constructions of Gray codes satisfying certain additional properties have been widely studied \cite{knuth2005art,Savage}.

\medskip
\textbf{1.1 Extending matching to Hamilton cycle.} Ruskey and Savage \cite{Ruskey} asked the~following question which is still open.

\begin{conjecture}[Ruskey and Savage \cite{Ruskey}]
	\label{conjecture:ruskey-savage}
	Every matching of the~hypercube~$Q_n$ can be extended into a~Hamilton cycle where $n \ge 2$.
\end{conjecture}

One natural step toward this problem is considering perfect matchings only, where each vertex is incident to exactly one edge.
Kreweras \cite{Kreweras} conjectured that every perfect matching in the~$n$-dimensional hypercube with $n \ge 2$ extends to a~Hamilton cycle.
This conjecture was popularized by Knuth \cite{knuth2005art} and proven by Fink \cite{kreweras1}.
The~proof of Kreweras' conjecture actually provides the following stronger statement where $K(Q_n)$ denote the~complete graph on vertices of~$Q_n$.

\begin{theorem}[Fink~\cite{kreweras1}]
	\label{theorem:perfect-matching}
	For $n \geq 2$ and every perfect matching~$M$ of $K(Q_n)$, there exists a~perfect matching $N$ of~$Q_n$, such that~$M \cup N$ is a~Hamilton cycle of $K(Q_n)$.
\end{theorem}

This result inspired several generalizations \cite{thomassen2015extending,gregor2009perfect}, e.g. the~authors of \cite{thomassen2015extending} showed that Kreweras' conjecture also holds for sparse spanning regular subgraphs of hypercubes.
Using the~following theorem, Alahmadi et al. \cite{alahmadi2015extending} proved that for every perfect matching~$M$ in~$Q_n$ the~set of edges $e$ of~$Q_n$ such that~$M \cup \set{e}$ cannot be extended into a~Hamilton cycle forms a~matching of~$Q_n$.

\begin{theorem}[Alahmadi et al.~\cite{alahmadi2015extending}, Theorem 4]
	\label{theorem:alahmadi}
	Let $n \geq 2$ be an~integer and let~$M$ be a~perfect matching on~$Q_n$.
	Let $e_1, e_2 \in E(Q_n) \setminus M$ incident with the~same vertex.
	Then~$Q_n$ has a~perfect matching $N$ such that $N$ contains one of $e_1, e_2$ and~$M \cup N$ is a~Hamilton cycle of~$Q_n$.
\end{theorem}

Shujia Wang and Fan Wang~\cite{wang24} presented another relaxation of the~Ruskey and Savage conjecture.
For vertices $u$ and $v$ connected by an~edge $u v$, we call the~position $i \in [n]$ where they differ the~\defi{direction} of that edge, with both \defi{endpoints} $u$ and $v$ being \defi{incident} to and \defi{covered} by this edge.
A~matching~$M$ in~$Q_n$ \defi{spans in direction $i$} if it contains at least one edge in that direction, and it \defi{spans at most~$d$ directions} if its edges lie in at most~$d$ distinct directions.

\begin{theorem}[Shujia Wang, Fan Wang~\cite{wang24}]
	\label{theorem:wang}
	For $n \geq 5$, let~$M$ be a~matching on the~hypercube~$Q_n$, such that $|M| < 10 \cdot 2^{n-5}$ and~$M$ spans at most~$5$ directions.
	Then~$M$ can be extended to a~Hamilton cycle.
\end{theorem}

We expand upon their idea, considering only matchings spanning at most~$5$ dimensions, and prove that the~size restriction on~$M$ is unnecessary.

\begin{theorem}
	\label{theorem:cycle-dim-5}
	For $n \geq 2$, let~$M$ be a matching on the~hypercube~$Q_n$ that spans at most~$5$ directions.
	Then~$M$ can be extended to a~Hamilton cycle.
\end{theorem}

Our proof of Theorem \ref{theorem:cycle-dim-5} works for arbitrary many directions $d$ assuming the~following conjecture holds.
Each vertex has a~\defi{parity} determined by the~parity of number of ones in its binary string representation.

\begin{conjecture}
	\label{conjecture:main-cycle}
	Given a~matching~$M$ on the~hypercube~$Q_n$ and two vertices $u, v \in V(Q_n)$ of opposite parity, where $u$ is not covered by~$M$, there exists a~Hamilton path extending~$M$ on~$Q_n$ and connecting $u$ and $v$ where $n \geq 2$.
\end{conjecture}

Observation \ref{observation:main-cycle-implies-ruskey-savage} proves that Conjecture \ref{conjecture:main-cycle} implies Conjecture \ref{conjecture:ruskey-savage}.
We verify this conjecture for $2 \le d \le 5$ in Observation \ref{observation:main-cycle-verification} using a computer program\footnote{\url{https://gitlab.com/jirka.fink/matching_of_five_directions}} which provides us Theorem \ref{theorem:cycle-dim-5}.
For general $d$, we have the~following conditional statement.

\begin{theorem}
	\label{theorem:hamilton-cycle}
	For $n \geq d \geq 2$, let~$M$ be a~matching on the~hypercube~$Q_n$ that spans at most~$d$ directions.
	If Conjecture~\ref{conjecture:main-cycle} holds for dimension $d$, then~$M$ can be extended to a~Hamilton cycle.
\end{theorem}

\medskip
\textbf{1.2 Extending matching to Hamilton path.}
It is well-known that the~hypercube~$Q_n$ contains a~Hamilton path between any two prescribed end vertices of opposite parity~\cite{simmons1977almost}.
Gregor, Novotn\'{y}, and {\v{S}}krekovski~\cite{gregor2017} considered the~problem of extending a~perfect matching of~$Q_n$ to a~Hamilton path between two prescribed end vertices with opposite parity.
Their proof works in the~more general setting of the~complete bipartite graph~$B(Q_n)$ which contains an~edge between all pairs of vertices of~$Q_n$ with opposite parities.
For a~matching~$M$ of~$Q_n$ and one of its vertices~$u$, we write $u^M$ for the~other end vertex of the~edge of~$M$ incident with~$u$.
A~\defi{layer} in direction $i \in [n]$ consists of all edges in direction $i$, and a~\defi{half-layer} is a~subset of a~layer containing only edges of specific parity.

\begin{theorem}[Gregor, Novotn\'{y}, and \v{S}krekovski~\cite{gregor2017}, Theorem 2]
	\label{theorem:skrekovski}
	Let $n\ge 5$, and let~$u, v$ be two vertices of opposite parity in~$Q_n$.
	A~perfect matching~$M$ of~$B(Q_n)$ with~$uv\notin M$ can be extended to a~Hamilton path, by edges from~$Q_n$, with end vertices~$u$ and~$v$ if and only if and $(M\setminus\{uu^M,vv^M\})\cup\{u^Mv^M\}$ contains no half-layers.
\end{theorem}

Motivated by this result, we ask when a~matching of~$Q_n$ can be extended into a~Hamilton path between prescribed vertices; and we conjecture the~following statement.
We call a~set of edges an~\defi{almost half-layer} if it is missing one edge to become a~half-layer.
For a~vertex~$u$, we denote by~$u^i$ its neighbor in direction $i$.
An~almost half-layer is \defi{$u$-avoiding} if it lacks only one edge incident to~$u$ to form a~complete half-layer.

\begin{conjecture}
	\label{conjecture:main-path}
	Let $n \geq 5$ and~$u, v \in V(Q_n)$ be vertices of opposite parity.
	For a~matching~$M$ on~$Q_n$, there exists a~Hamilton path between~$u$ and $v$ extending~$M$ if and only if none of these conditions holds:

	\begin{enumerate}[label={C\theenumi.}]
		\item $M$ contains a~half-layer and~$u$ and $v$ are covered by it,
		\item $M$ contains a~$u$-avoiding almost half-layer in the~direction $i$ where $v = u^i$, $j \neq i$ and \mbox{$u u^j, v v^j \in M$},
		\item $u v \in M$.
	\end{enumerate}
\end{conjecture}

\begin{figure}
	\begin{center}
		\includegraphics[width=0.5\textwidth]{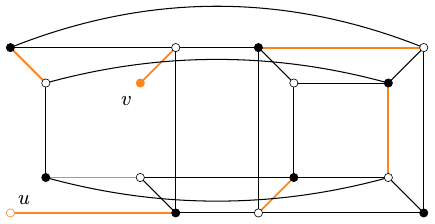}
	\end{center}
	\caption{
		The~forbidden case for $n=4$.
		Orange edges are the~matching~$M$ and the~orange vertices are the~vertices $u$, $v$ which cannot be together extended to a~Hamilton path extending~$M$ between $u$ and $v$.
		Black edges are the~remaining usable edges in~$Q_4$.
	}
	\label{figure:forbidden-matchings}
\end{figure}

Observation~\ref{observation:c-conditions-necessity} proves that all C-conditions are necessary for the~existence of a~Hamilton path, and our main conjecture is that they are also sufficient.
Note that the~requirement $n \ge 5$ is necessary, since counterexamples for $n=4$ presented in \cite{gregor2017} applies here as well, with one extra counterexample shown in Figure~\ref{figure:forbidden-matchings}.
Furthermore, Dvo\v{r}\'ak and Fink \cite{dvorak2016} constructed a~matching of $B(Q_n)$ of size $\Theta(2^n/\sqrt{n})$ which cannot be extended into a~Hamilton path for $n \ge 9$.
Using a similar construction, Observation \ref{observation:BQn-not-extendable} presents a~matching of $B(Q_n)$ which cannot be extended to a~Hamilton cycle between any two vertices for $n \ge 9$, so we restrict Conjecture \ref{conjecture:main-path} to consider only matchings of~$Q_n$.
Observation \ref{observation:conjecture-path-implies-cycle} proves that Conjecture \ref{conjecture:main-path} implies \ref{conjecture:main-cycle}.

We consider the~following simplification of Conjecture \ref{conjecture:main-path} for maximal matchings, where a~matching is maximal if it cannot be extended by any edge.

\begin{conjecture}
	\label{conjecture:maximal}
	Let $n\ge 5$, and let~$u, v$ be two vertices of opposite parity in~$Q_n$.
	A~maximal matching~$M$ of~$Q_n$ with~$uv\notin M$ can be extended to a~Hamilton path with end vertices~$u$ and~$v$ if and only if and $(M\setminus\{uu^M,vv^M\})\cup\{u^Mv^M\}$ contains no half-layers.
\end{conjecture}

Observation \ref{observation:maximal-implies} proves that Conjecture \ref{conjecture:main-path} implies Conjecture~\ref{conjecture:maximal}.
This shows the~relation between Conjecture \ref{conjecture:main-path} and Theorem \ref{theorem:skrekovski}, and it implies that Conjecture \ref{conjecture:main-path} holds for perfect matchings and matchings extendable to a~perfect one.

To further support the~main Conjecture \ref{conjecture:main-path}, we verify it for matching spanning at most $d=5$ directions in Observation \ref{observation:main-path-verification}.
Similarly as in Theorem \ref{theorem:hamilton-cycle}, we prove the~following conditional version of Conjecture \ref{conjecture:main-path} for matching spanning at most $d$ directions, and verify the~condition for $d=5$.

\begin{theorem}
	\label{theorem:hamilton-path}
	Let $n \geq d \geq 5$, and~$u, v \in Q_n$ be vertices of opposite parity and~$M$ be a~matching on~$Q_n$ spanning at most~$d$ directions.
	Assuming Conjecture~\ref{conjecture:main-path-two} holds for dimension $d$, then~$M$ can be extended into a~Hamilton path between~$u$ and $v$ if and only if none of the~C-conditions hold for~$M$.
\end{theorem}

The~extra condition in Theorem \ref{theorem:hamilton-path} is the~following and Conjecture~\ref{conjecture:main-path} is a~direct consequence of it.

\begin{conjecture}
	\label{conjecture:main-path-two}
	Let $n \geq 5$, and let~$u, v \in Q_n$ be vertices of opposite parity and~$M$ be a~matching on~$Q_n$.
	Then if and only if none of the~C-conditions hold for~$M$ there exists two distinct Hamilton paths between~$u$ and $v$ extending~$M$.
\end{conjecture}

Observation \ref{observation:main-path-verification} verifies that Conjecture \ref{conjecture:main-path-two} holds for $d = 5$ which implies the~following statement.

\begin{theorem}
	\label{theorem:hamilton-path-5}
	Let $n \geq 5$, and~$u, v \in Q_n$ be vertices of opposite parity and~$M$ be a~matching on~$Q_n$ spanning at most~$5$ directions.
	Then $M$ can be extended into a~Hamilton path between~$u$ and $v$ if and only if none of the~C-conditions hold for~$M$.
\end{theorem}

\medskip
\textbf{1.3 Related results on extending edges to cycles.}
The~research on Hamilton cycles in hypercubes satisfying certain additional properties has received considerable attention; see e.g. the~survey by Savage \cite{Savage}.
They are motivated by applications in computer networks, where the~hypercube is frequently used as a~network topology with a~number of desirable properties, such as small degree and diameter.

Vandenbussche and West \cite{vandenbussche2013extensions} who conjectured that every matching can be extended into \mbox{a~2-factor}, and the~conjecture was proven by Fink \cite{2factor}.
Fink and M\"utze \cite{fink2024matchings} proved that every matching of~$Q_n$ can be extended into a~cycle that visits at least a~2/3-fraction of all vertices.

Dvo\v{r}\'{a}k \cite{Dvorak_cycle_prescribed} showed that any set of at most $2n-3$ edges of~$Q_n$, $n \ge 2$, that induces vertex-disjoint paths is contained in a~Hamilton cycle.
Dvo\v{r}ák and Gregor \cite{dvovrak2007hamiltonian} proved that for any set $E$ of at most $2d-4$ edges in~$Q_n$, $n \ge 5$, that form disjoint paths and two vertices x and y of opposite parity that are neither internal vertices of the~paths nor end vertices of the~same path, there is a~Hamilton path with end vertices $x$ and $y$ that contains all of $E$.

Dimitrov et al. \cite{skrekovski2009gray} proved that the~hypercube~$Q_n$ contains a~Hamilton cycle avoiding a~given matching except a~forbidden configuration.

Other related results are nicely surveyed in the introduction of recently published paper \cite{fink2024matchings}.

\medskip
\textbf{1.4 Overview of our approach.}
In this paper, we consider matchings~$M$ in~$Q_n$ spanning at most $d \le n$ directions.
So, it is natural to split~$Q_n$ into $2^{n-d}$ subcubes of dimension $d$ so that every edge of~$M$ belongs into one of these subcubes.
We find a~Hamilton cycle extending~$M$ by finding cycles extending~$M$ in every subcube so that these cycles can be interconnected into a~single Hamilton cycle.
A~similar approach was presented by Gregor \cite{gregor2009perfect} who proved that every perfect matching $N$ of~$Q_n$ can be extended into a~Hamiltonian cycle by adding only edges $I$ inside prescribed subcubes if and only if $N \cup I$ is connected.
The~connectivity condition is clearly necessary, so in our case where~$M \subseteq I$, we have to use some edges between subcubes in the~Hamilton cycle extending~$M$, but our approach use only few of them.

%In section \ref{sec:relation} we show some important relations between conjectures that were discussed in introduction.
%Extending matchings into a~cycle and a~path is presented in Section \ref{sec:cycle} and \ref{sec:path}, resp.
%Programs verifying our statements for $d=5$ is described in Section \ref{sec:computer}.\todo{Probably need to mention section \ref{sec:matching-not-extendable}. Does it make a~difference that it is in appendix?}

Extending matchings into a~cycle and a~path is presented in Section \ref{sec:cycle} and \ref{sec:path}, resp.
In Appendix \ref{sec:relation} we show some important relations between conjectures that were discussed in introduction.
Programs verifying our statements for $d=5$ is described in Appendix \ref{sec:computer}.

\section{Extending to a~Hamilton cycle}
\label{sec:cycle}

In this section, we show how to extend a~matching~$M$ on the~hypercube~$Q_n$ that spans at most~$d$ directions into a~Hamilton cycle.
First, we prove Theorem \ref{theorem:hamilton-cycle} giving the~extension for $n \geq d \geq 2$, but it assumes that Conjecture~\ref{conjecture:main-cycle} holds for dimension $d$.
Then, using Observation \ref{observation:main-cycle-verification} which verifies Conjecture~\ref{conjecture:main-cycle} for $d=5$, we prove Theorem \ref{theorem:cycle-dim-5}.

First we introduce some notation.
The~\defi{Cartesian product} $G \cp H$ of graphs $G$ and $H$ has vertex set $V(G) \times V(H)$, with edges between vertices~$(u, v)$ and~$(u', v')$ if either~$u = u'$ and~$vv' \in E(H)$, or~$v = v'$ and~$uu' \in E(G)$.
Note that $Q_n = Q_{n-d} \cp Q_d$.
We call the~copies of $H$ in $G$ \defi{canonical copies}.
In this paper, we consider matchings on~$Q_n$ that span at most~$d$ directions, and we can assume without loss of generality that these directions are the~last $d$ directions.
Therefore, every edge $uv$ of these matchings belongs to one canonical copy of~$Q_d$ in~$Q_n$, and this canonical copy is uniquely identified by the~first $n-d$ directions $u$ which are the~same as in $v$.
For a~vertex $u \in V(Q_n)$, we denote by $u^L$ and $u^R$ its first (Left) $n-d$ and last (Right) $d$ coordinates, respectively.
Let~$Q_n[u^L]$ be subgraph of~$Q_n$ induced by vertices $w$ such that $w^L = u^L$.
In another works, $Q_n[u^L]$ is the~canonical copy of~$Q_n$ containing $u$.
Similarly for a~matching~$M$ spanning only the~last $d$ directions, $M[u^L]$ are edges of~$M$ in~$Q_n[u^L]$.

\begin{figure}
	\begin{center}
		\includegraphics[width=0.5\textwidth]{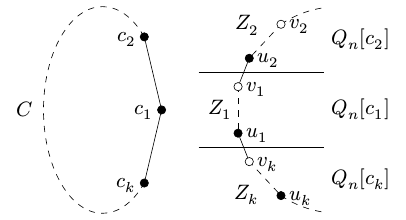}
	\end{center}
	\caption{
		Illustration of proof of Theorem~\ref{theorem:hamilton-cycle}.
		On the~left is a~Hamilton cycle $C = (c_1, \dots, c_k)$ on~$Q_{n-d}$ and on the~right the~solid lines are edges between subcubes and the~dashed lines are the~Hamilton paths on the~subcubes.
	}
	\label{figure:hamilton-cycle-proof}
\end{figure}

\begin{proof}[Proof of Theorem~\ref{theorem:hamilton-cycle}]
	If $n = d$, the~theorem is a~direct consequence of Observation~\ref{observation:main-cycle-implies-ruskey-savage}.
	Assuming $n > d$, we can without loss of generality consider only matchings that span the~last~$d$ directions.
	If~$M$ is perfect, we invoke Theorem~\ref{theorem:perfect-matching} for immediate proof.
	Otherwise, for a~non-perfect matching, there is vertex $v_1$ avoided by~$M$.
	Let $k = 2^{n-d}$, and consider a~Hamilton cycle $C = (c_1, c_2, \dots, c_k)$ on~$Q_{n-d}$, such that~$Q_n[c_1]$ is the~subcube containing $v_1$.

	The~objective is for each $i \in [k]$ to construct a~Hamilton path $Z_i$ in~$Q_n[c_i]$ extending~$M[c_i]$ and connecting some vertices $u_i$ and $v_i$.
	Also ensuring $u_i$ and $v_{i-1}$ are adjacent for $i \in \set{2, 3, \dots k}$, as well as $v_k$ to $u_1$.
	These paths, together with edges $u_i v_{i-1}$ for $i \in \set{2, 3, \dots k}$ and $v_k u_1$, complete a~Hamilton cycle extending~$M$ which implies the~statement.
	Figure~\ref{figure:hamilton-cycle-proof} shows how the~proof works.

	Paths $Z_i$ are constructed inductively, starting from vertex $v_1$ avoided by~$M$ in~$Q_n[c_1]$.
	For each $i \in \set{2, 3, \dots, k}$, $u_i$ is chosen as $v_{i-1}$ sole neighbor in~$Q_n[c_i]$.
	From Observation~\ref{observation:main-cycle-implies-ruskey-savage} each~$Q_n[c_i]$ contains a~Hamilton cycle $C_i$ extending the~matching~$M[c_i]$.
	Selecting $v_i$ among $u_i$'s neighbors in $C_i$, such that $u_i v_i \notin M$, creates a~Hamilton path $Z_i = E(C_i) \setminus \set{ u_i v_i }$ between $u_i$ and $v_i$ extending~$M[c_i]$.

	Upon reaching $v_k$, $u_1$ is identified as the~sole $v_k$ neighbor in~$Q_n[c_1]$.
	Applying Conjecture~\ref{conjecture:main-cycle} to~$M[c_1]$ and vertices $u_1$ and $v_1$ yields a~Hamilton path $Z_1$ from $u_1$ to $v_1$, extending~$M[c_1]$ and completing the~proof.
\end{proof}

\begin{proof}[Proof of Theorem~\ref{theorem:cycle-dim-5}]
	The~theorem directly follows from Observation \ref{observation:main-cycle-verification} and Theorem \ref{theorem:hamilton-cycle}.
\end{proof}

\section{Extending to Hamilton path}
\label{sec:path}

The~extension of matchings to Hamilton paths presents additional challenges beyond the~cycle case.
Half-layers create rigid structures that constrain possible paths through the~hypercube, making them natural obstacles to path extension.
This phenomenon was observed by Gregor, Novotn\'{y} and \v{S}krekovski~\cite{gregor2017} in their study of path extensions in bipartite graphs derived from hypercubes.

Their work, presented in Theorem~\ref{theorem:skrekovski}, characterizes path extensions from perfect matchings in $B(Q_n)$.
They also showed an~equivalent result.
For a~graph $G$, let $G - X$ denote the~graph obtained by removing vertices in set $X$ and all incident edges.

\begin{theorem}[Gregor, Novotn\'{y}, and \v{S}krekovski~\cite{gregor2017}, Theorem 3]
	\label{theorem:skrekovski-2}
	Let $u, v$ be vertices of opposite parity in~$Q_n, n \geq5$, and let~$M$ be a~perfect matching of $B(Q_n - \{u, v\})$.
	There exists a~matching $N \subseteq E(Q_n)$ such that $N \cup M$ forms a~Hamilton cycle of $B(Q_n - \{u, v\})$ if and only if~$M$ does not contain a~half-layer.
\end{theorem}

The~dimension requirement $n \geq 5$ proves essential, providing sufficient flexibility for routing Hamilton cycles through the~hypercube while accommodating matching constraints.
Smaller dimensions lack this flexibility, as demonstrated by $n = 4$ where certain matching configurations provably~\cite{gregor2017} prevent Hamilton path extensions.

This characterization builds on a~fundamental structural property of half-layers:

\begin{lemma}[Gregor, Novotn\'{y}, and \v{S}krekovski~\cite{gregor2017}, Lemma 4]
	\label{lemma:skrekovski}
	Let $H_i$ and $H_j$ be half-layers of~$Q_n$ of different directions $i$ and $j$, respectively.
	Then, every edge of $H_i$ is incident with some edge of $H_j$.
\end{lemma}

This intersection property has important implications: if a~matching contains a~half-layer in some direction $i$, it cannot contain even a~single edge of a~half-layer in any other direction $j$, as these edges would necessarily share endpoints.
We can see that if a~perfect matching contains two half-layers, they both must be of the~same direction and of opposite parity.
We also need the~following lemma regarding almost half-layers:

\begin{lemma}[Gregor, Novotn\'{y}, and \v{S}krekovski~\cite{gregor2017}, Lemma 6]
	\label{lemma:skrekovski-almost}
	Let~$M$ be a~matching of $K(Q_n), n \geq 4$. Then,~$M$ contains almost half-layer in at most one direction.
\end{lemma}

The~pathological nature of half-layers is most evident when considering the~C-conditions:

\begin{figure}
	\begin{center}
		\includegraphics[width=0.4\textwidth]{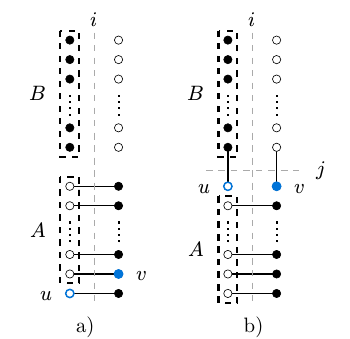}
	\end{center}
	\caption{Necessary conditions for the~existence of Hamilton paths in a~hypercube, a) condition C1, $u$ and $v$ covered by a~half-layer b) condition C2,~$M$ contains $u$ and $v$ avoiding almost half-layer in direction $i$ and $u u^j, v v^j \in M, j \neq i$.}
	\label{figure:hypercube-path-conjecture-conditions}
\end{figure}

\begin{observation}
	\label{observation:c-conditions-necessity}
	If any C-condition from Conjecture~\ref{conjecture:main-path} holds, no Hamilton path between~$u$ and~$v$ extends a~matching~$M$.
\end{observation}

\begin{proof}
	Negation of condition C3 is trivially necessary.
	We prove only the~necessity of case C2, since the~C1's proof is similar.
	Both cases are illustrated in Figure~\ref{figure:hypercube-path-conjecture-conditions}.

	For C2, let $i,j$ be as defined as in the~condition, with $H$ denoting the~$u$-avoiding almost half-layer in direction $i$.
	Without loss of generality, assume~$u$ is even, and~$v$ is odd.
	Let $A$ be the~set of even vertices, except~$u$, covered by~$H$ and~$B$ the~set of odd vertices not covered by~$H$, except~$v$.
	For path to exist between~$u$ and~$v$ extending~$M$, each vertex in $B$ must have an~edge incident to it, that is not in the~direction of $i$.
	That holds also for $u^j$, since the~edge $u^j v^j$ would close the~path prematurely.
	So each vertex in $B$ needs its counterpart from~$A$.
	However, $|B| = 2^{n-2}$ while $|A| = 2^{n-2}-1$, making such a~path construction impossible.
\end{proof}

The~conditions C1 and C2 have structural implications: when either holds, $M$ must contain either a~half-layer or an~almost half-layer in some direction~$i$ and from Lemma~\ref{lemma:skrekovski} and Lemma~\ref{lemma:skrekovski-almost} only in that direction.
Thus, only C1 or C2 may hold at once.
Moreover, these conditions naturally suggest how~$M$ could be completed to a~perfect matching: any extension must use edges in direction~$i$ to maintain the~half-layer or the~almost half-layer structure.
We call this completion the~\defi{unique extension of~$M$}, as it aligns with the~conditions characterized in Theorem~\ref{theorem:skrekovski}.
Since any Hamilton path extending the~unique extension must use all its edges, such a~path also extends~$M$.

The~following lemma characterizes when there exists a~Hamilton path between a~vertex and one of its neighbors:

\begin{lemma}
	\label{lemma:three-neighbors}
	Let $d \geq 5$, $u \in V(Q_d)$, and~$M$ be a~matching of~$Q_d$.
	Assuming Conjecture~\ref{conjecture:main-path} holds for dimension $d$, then for at most two of the~neighbors of $u$ there does not exist a~Hamilton path extending~$M$ connecting them to $u$.
\end{lemma}
\begin{proof}
	Consider which neighbors of $u$ are forbidden by C-conditions.
	Under condition C1, $u$ must be covered by a~half-layer in direction~$i$.
	This forbids only neighbor $u^i$ both by C1 and C3 because~$uu^i \in M$ lies in the~half-layer, while all other neighbors remain valid as they lie on the~correct side of the~half-layer.

	Under condition C2 with neighbor $v = u^i$, $M$ contains a~$u$-avoiding almost half-layer in direction $i$ and edges $uu^j, vv^j \in M$ for some $j \neq i$.
	This forbids exactly two neighbors: $v = u^i$ by C2 and $u^j$ by C3, since $uu^j \in M$.

	As C1 and C2 cannot hold simultaneously, at most two neighbors are forbidden, implying any set of three neighbors contains at least one admitting a~Hamilton path extending~$M$.
\end{proof}

Beyond examining paths to neighbors, we consider paths to all vertices of opposite parity from a~given vertex~$u$.
Let $T^M_u$ denote the~set of vertices~$v$ where a~Hamilton path exists between~$u$ and~$v$ extending~$M$.
This notation enables us to quantify exactly how many endpoints are available for Hamilton paths:

\begin{observation}
	\label{observation:available-endpoints}
	Let $n \geq 5$,~$M$ be a~matching of~$Q_n$, and $u \in V(Q_n)$.
	Assuming Conjecture~\ref{conjecture:main-path} holds for $n$, then $|T^M_u| \geq 2^{n-2}$.
	Moreover, $|T^M_u| = 2^{n-2}$ if and only if~$M$ contains a~half-layer $H$ covering $u$, in which case $T^M_u$ consists exactly of vertices with parity opposite to $u$ not covered by $H$.
	Otherwise $|T^M_u| \geq 2^{n - 1} - 2$.
\end{observation}

\begin{proof}
	If $u$ is not covered by~$M$, none of the~C-conditions can hold between $u$ and any vertex of opposite parity.
	Thus $|T^M_u| = 2^{n-1}$.
	For the~remainder of the~proof, we assume $u$ is covered by~$M$ and analyze three cases.

	First, suppose~$M$ contains a~half-layer $H$ covering $u$.
	By condition C1, $T^M_u$ consists exactly of vertices of opposite parity to $u$ that are not covered by $H$, giving $|T^M_u| = 2^{n-2}$.

	Next, suppose~$M$ contains a~$u$-avoiding almost half-layer in direction $i$.
	Then $T^M_u$ excludes exactly two vertices of opposite parity.
	$u^M$, that is forbidden by C3 since $u u^M \in M$ and $u^i$, that is forbidden by C2.
	Therefore $|T^M_u| = 2^{n-1} - 2$.

	Finally, if neither C1 nor C2 holds, only $u^M$ is forbidden by C3, yielding $|T^M_u| = 2^{n-1} - 1$.

	Therefore, $|T^M_u| = 2^{n-2}$ occurs if and only if C1 holds.
	In this case, by Lemma~\ref{lemma:skrekovski}, the~direction and parity of $H$ are uniquely determined by the~complement of $T^M_u$.
	In all other cases, $|T^M_u| \geq 2^{n-1} - 2$ for $n \geq 5$.
\end{proof}

This observation gives us an~important information about the~Hamilton paths extending matching~$M$ and stemming from a~specific vertex $u$.
Either we have many available endpoints ($|T^M_u| \geq 2^{n-1} - 2$) or $M$ contains a~specific half-layer covering $u$ that determines the~endpoint set as all of the~vertices of opposite parity to $u$ not covered by the~said half-layer.

To analyze path constructions systematically, we introduce cube sequences as our primary structural tool.
A~\defi{sequence of $m$ cubes of dimension $d$} is the~graph $P_m \cp Q_d$, where $P_m$ is the~path on $m$ vertices.
For a~path $Y = (y_1, \dots, y_m)$ in~$Q_{n-d}$ where $m \leq 2^{n-d}$, its induced sequence in~$Q_n$ consists of subcubes $G_i = Q_n[y_i]$ and edges between them, with corresponding matchings~$M_i = M[y_i]$ for $i \in [m]$.
We refer to this as the~\defi{induced sequence of cubes and matchings of path $Y$}.

Any cube sequence embeds as a~subgraph into~$Q_n$ for sufficiently large $n$, allowing us to extend our hypercube notation to these sequences.
Vertices and edges in cube sequences exhibit a~dual nature: they possess both local properties within their canonical copy of~$Q_d$ and global properties from their embedding in~$Q_n$.
For instance, in $Q_n = Q_{n-d} \cp Q_d$, a~vertex $u \in V(Q_n)$ with $u^L = \bm{1} \in V(Q_{n-d})$ and $u^R = \bm{0} \in V(Q_d)$ is locally even in its~$Q_d$ copy due to $u^R = \bm{0}$, while its global parity in~$Q_n$ is determined by $u^L = \bm{1}$.

A~\defi{sequence of half-layers in $P_m \cp Q_d$} consists of half-layers in each canonical copy of~$Q_d$, all sharing the~same direction and global parity.
These sequences exhibit an~important alternation property: when traversing consecutive cubes in the~sequence of cubes, half-layers of the~same global parity must alternate between odd and even local parity within their respective~$Q_d$ copies—a~property inherent to the~embedding structure.

For a~subset $W \subseteq V(Q_n)$, we define its projection onto the~last $d$ coordinates as $W^R = \set{w^R \mid w \in W}$, where $W^R \subseteq V(Q_d)$.
We use this projection notation to distinguish between global vertices in~$Q_n$ and their local representations in canonical copies of~$Q_d$.

We are ready to prove our first key structural result:

\begin{lemma}
	\label{lemma:two-cube-path}
	Let $d \geq 5$ and $G = P_2 \cp Q_d$ with canonical copies $G_1, G_2$ of~$Q_d$ containing matchings~$M_1, M_2$ respectively.
	Let $u \in G_1, v \in G_2$ have different global parities.
	Assuming Conjecture~\ref{conjecture:main-path} holds for dimension $d$, there exists a~Hamilton path between $u$ and $v$ extending~$M = M_1 \cup M_2$ if and only if no sequence of half-layers in $G$ contained in~$M$ covers both $u$ and $v$.
	Moreover, when such a~path exists, there are at least $2^{d-3}$ distinct Hamilton paths, determined by different edges between $G_1$ and $G_2$ satisfying these conditions.
\end{lemma}

\begin{proof}
	We begin by proving necessity.
	If $G$ contains a~sequence of half-layers in~$M$ covering both~$u$ and~$v$, this sequence corresponds to a~half-layer in~$Q_{d+1}$ covering~$u$ and~$v$.
	By Observation~\ref{observation:c-conditions-necessity}, no Hamilton path can exist.

	For sufficiency, assume no sequence of half-layers in~$M$ covers both~$u$ and~$v$.
	Let~$T^{M_1}_u$ and~$T^{M_2}_v$ be their respective sets of available endpoints in their respective copies $G_1$ and $G_2$.
	Given $w \in (T^{M_1}_u)^R \cap (T^{M_2}_v)^R$, let $w_1 \in V(G_1), w_2 \in V(G_2)$ such that $w_1^R = w$ and $w_2^R = w$.
	By definition of~$T^{M_1}_u$ and~$T^{M_2}_v$, there exist Hamilton paths $Z_1$ from $u$ to $w_1$ extending~$M_1$ in $G_1$ and $Z_2$ from $v$ to~$w_2$ extending~$M_2$ in~$G_2$.
	The~path $Z_1 \cup \set{w_1w_2} \cup Z_2$ forms a~Hamilton path in $G$ extending~$M$.

	To establish the~bound on the~number of distinct Hamilton paths, we consider several cases.
	When $u$ or $v$ is not covered by a~half-layer in~$M$, say $u$, then by Observation~\ref{observation:available-endpoints}, $|T^{M_1}_u| \geq 2^{d-1} - 2$.
	Therefore $|(T^{M_1}_u)^R \cap (T^{M_2}_v)^R| \geq |T^{M_2}_v| - 2 \geq 2^{d-2} - 2 \geq 2^{d-3}$ for $d \geq 5$.

	When both vertices are covered by half-layers, these must differ either in direction or global parity.
	In the~case of same direction but different global parity, $(T^{M_1}_u)^R = (T^{M_2}_v)^R$ with $|T^{M_1}_u| = 2^{d-2}$.
	When the~directions differ, the~overlap has size exactly $2^{d-3}$, as this corresponds to fixing two coordinates for vertices of specific parity.
	Thus in all cases, we obtain at least $2^{d-3}$ distinct choices for $w$, each yielding a~distinct Hamilton path.
\end{proof}

We now generalize Lemma~\ref{lemma:two-cube-path} for longer sequences of hypercubes:

\begin{lemma}
	\label{lemma:path-sequence}
	Let $d \geq 5$ and $G = P_m \cp Q_d$ with $m \geq 2$.
	For $i \in [m]$, let~$M_i$ be a~matching of the~$i$-th canonical copy $G_i$ of~$Q_d$, and let~$M = \cup_{i \in [m]} M_i$.
	For vertices $u \in G_1, v \in G_m$ of opposite global parity, there exists a~Hamilton path between $u$ and $v$ extending~$M$ if and only if no sequence of half-layers in $G$ contained in~$M$ covers both $u$ and $v$, assuming Conjecture~\ref{conjecture:main-path} holds for dimension $d$.
	Moreover, when connecting consecutive canonical copies of~$Q_d$, we can choose from $2^{d-3}$ different edges.
\end{lemma}

\begin{proof}
	We proceed by induction on $m$.
	The~base case $m=2$ follows from Lemma~\ref{lemma:two-cube-path}.
	For the~inductive step, let $m > 2$ and assume the~statement holds for $m-1$.
	Let $G' = P_{m-1} \cp Q_d$ and~$M'$ be the~restriction of~$M$ to $G'$.
	We denote by $T^{M'}_u$ the~set of vertices in $G_{m-1}$ that are endpoints of Hamilton paths in $G'$ extending~$M'$ and starting at $u$.
	Similarly, let $T^{M_m}_v$ be the~set of vertices in $G_m$ that are endpoints of Hamilton paths extending~$M_m$ and starting at $v$.

	Notice if $u$ is covered by a~sequence of half-layers in~$M'$, $T^{M'}_u$ consists of vertices of opposite global parity to $u$ not covered by that sequence.
	Thus if there exists a~sequence of half-layers in~$M$ covering both $u$ and $v$, then $(T^{M'}_u)^R$ and $(T^{M_m}_v)^R$ are disjoint since vertices in $(T^{M'}_u)^R$ correspond to vertices covered by the~half-layer covering $v$ in $G_m$, thus no Hamilton path exists.

	If no sequence of half-layers covers $u$ in $G'$, then $T^{M'}_u$ contains all vertices of opposite global parity to $u$ in $G_{m-1}$.
	Since $|T^{M_m}_v| \geq 2^{d-2}$ by Observation~\ref{observation:available-endpoints}, $|(T^{M'}_u)^R \cap (T^{M_m}_v)^R| \geq 2^{d-2}$, and a~Hamilton path exists by the~construction in Lemma~\ref{lemma:two-cube-path}.

	When a~sequence of half-layers covers~$u$ in~$G'$ but no sequence of half-layers in~$G$ covers both~$u$ and~$v$, we have $|T^{M'}_u| = 2^{d-2}$.
	For $|T^{M_m}_v| = 2^{d-2}$, the~half-layer in~$M_m$ covering $v$ is not part of the~sequence covering $u$.
	If it is part of a~half-layer of the~same direction as the~sequence we get, $|(T^{M'}_u)^R \cap (T^{M_m}_v)^R| = 2^{d-2}$.
	If the~half-layer covering $v$ is of different direction than of the~sequence we get $|(T^{M'}_u)^R \cap (T^{M_m}_v)^R| = 2^{d-3}$, because the~intersection consists of vertices of specific parity with two positions fixed.
	When $|T^{M_m}_v| > 2^{d-2}$, then it must hold that $|T^{M_m}_v| \geq 2^{d-1} - 2$ and for the~intersection it holds that $|(T^{M'}_u)^R \cap (T^{M_m}_v)^R| \geq 2^{d-2} - 2$.
	In both cases, a~Hamilton path exists and we have at least $2^{d-3}$ of connections to choose when connecting the~$G'$ with $G_{m-1}$.
\end{proof}

Consider attaching an~additional cube of dimension $d$ to a~sequence of cubes of the~same dimension through which we have already threaded a~Hamilton path extending a~matching~$M$.
The~following lemma establishes the~conditions under which such an~attachment is possible for a~pair of cubes:

\begin{lemma}
	\label{lemma:cube-attachment}
	Let $d \geq 5$ and $G = P_2 \cp Q_d$ with canonical copies $G_1, G_2$ of~$Q_d$ containing matchings~$M_1, M_2$ respectively.
	Let~$Z_1$ and~$Z_2$ be two distinct Hamilton paths extending~$M_1$ in~$G_1$, where~$Z_1$ is between~$u$ and~$v_1$ and~$Z_2$ is between~$u$ and~$v_2$.
	Note that~$v_1$ and~$v_2$ are not necessarily distinct.
	Assuming Conjecture~\ref{conjecture:main-path} holds for dimension~$d$, there exists a~Hamilton path in~$G$ extending~$M_1 \cup M_2$ from~$u$ to at least one of~$\set{v_1, v_2}$.
\end{lemma}

\begin{proof}
	An~edge of~$M_1$ is called \defi{fixed}, and an~edge of $(Z_{v_1} \cup Z_{v_2}) \setminus M_1$ is said to be \defi{free}.
	For a~vertex~$a \in G_1$ let $a' \in G_2$ be its unique corresponding counterpart, such that $a^R = (a')^R$.
	For a~free edge $a b$, we say it \defi{admits a~Hamilton path} if there exists a~Hamilton path in $G_2$ extending~$M_2$ between the~corresponding vertices $a'$ and $b'$, which we denote by $Z'_{a'b'}$.
	Given such an~admitting edge $ab$ in some path~$Z_v$, we can construct a~Hamilton path in $G$ as $Z_v \setminus \set{ab} \cup \set{aa', bb'} \cup Z'_{a'b'}$.

	We claim that there exists at least one free edge admitting a~Hamilton path.
	Suppose toward contradiction that no such admitting free edge exists.

	Let $x$ be the~first vertex where paths $Z_{v_1}$ and $Z_{v_2}$ diverge, with $Z_{v_1}$ following edge $xy$ and $Z_{v_2}$ following $xz$.
	Note that neither $y$ nor $z$ can be an~end vertex $v_1$ or $v_2$, as this would force the~paths to be identical, since the~only vertex left to visit is $v_1 = v_2$.

	It is clear that both of the~vertices $y$ and $z$ must be incident to at least three edges from these paths, since they differ at these vertices.
	By Lemma~\ref{lemma:three-neighbors} at most two of them can be free, since if there would be three, at least one of them would admit a~Hamilton path.
	Since only one incident edge may be fixed we get two free edges and one fixed for both of these vertices.

	This structure combined with the~C-conditions, implies that for all of the~vertices~$x, y$ and~$z$ the~C2 and C3 conditions must hold and thus~$M_2$ contains two almost half-layers in some direction~$i$, avoiding the~vertices $x'$, $y'$ and $z'$.
	Let $w'$ be the~fourth vertex avoided by these almost half-layers, with $w$ its corresponding vertex in $G_1$.
	The~two other free edges coming from $y$ and $z$ are exactly $y w \in Z_{v_2}$ and $z w \in Z_{v_1}$.

	All remaining free edges must lie in direction $i$, since they do not admit a~Hamilton path, and thus alternate with fixed edges.
	This is forcing~$Z_{v_1}$ and~$Z_{v_2}$ to be identical except for their behavior at aforementioned vertices.
	The~path up to $x$ is the~same for both of the~paths, but the~vertices~$y, z$ and~$w$ must be connected somehow.
	If there would be a~path, common to both of~$Z_{v_1}$ and~$Z_{v_2}$, connecting $y$ and $w$ directly, without using $z$, there would be a~cycle together with the~edge~$y w$ in~$Z_{v_1}$.
	The~same applies to $z$ and $w$.
	Thus there must be a~path between $y$ and $z$ without using $w$.
	However, this creates a~path between vertices of the~same parity, using an~even number of edges alternating between free and fixed - an~impossibility since the~path starts and ends with a~fixed edge.

	This contradiction establishes that some free edge must admit a~Hamilton path, allowing us to construct the~desired path in $G$.
\end{proof}

Before proving Theorem~\ref{theorem:hamilton-path}, we must address the~case where $u$ and $v$ lie in the~same subcube and satisfy one of the~C-conditions, specifically either the~condition C1 or C2, since if they would satisfy the~condition C3 no Hamilton path would exist even in the~larger cube.
We use the~following result about perfect matchings in hypercubes:

\begin{lemma}
	\label{lemma:disjoint-paths}
	Let $n \geq 5$ and~$M$ be a~matching of~$Q_n$.
	Let $u, v \in V(Q_n)$ be vertices of different parity with $u v \notin M$.
	If condition C1 or C2 holds for $u$, $v$ and~$M$, then there exist vertices $z_u$ and $z_v$ of different parity such that $z_u$ is adjacent to $u$, $z_v$ is adjacent to $v$, but $z_u$ and $z_v$ are not adjacent to each other.
	Furthermore, there exist vertex-disjoint paths~$Z_u$ and~$Z_v$ where~$Z_u$ connects~$u$ to~$z_u$, $Z_v$ connects~$v$ to~$z_v$, and $Z_u \cup Z_v$ covers all vertices of~$Q_n$ while extending the~matching~$M$.
\end{lemma}

\begin{proof}[Proof of Lemma~\ref{lemma:disjoint-paths}]
	In this proof consider~$M$ being its unique extension to a~perfect matching.
	We proceed by analyzing each condition separately.

	First suppose condition C1 holds, so~$M$ contains a~layer in direction $i$.
	Since $u v \notin M$ and they are covered by the~same half-layer, their distance must exceed three.
	Therefore, there exists a~$k$ such that $u_k \neq v_k$, thus $u$ and $v$ lie in opposite subcubes, induced by the~$k$-th position.
	Let~$Q^k_0$ denote the~$(n-1)$-dimensional subcube induced by vertices with $0$ in position $k$, and~$Q^k_1$ the~subcube with 1 in position $k$.
	Let~$M^k_0$ and~$M^k_1$ be the~restrictions of~$M$ to these subcubes.

	In~$Q^k_0$ containing $u$, we apply Theorem~\ref{theorem:perfect-matching} to obtain a~cycle $C_u = (u, \dots, z_u)$, which forms our path $Z_u$.
	For~$Q^k_1$ containing $v$, we must ensure our path ends at a~vertex $z_v$ not adjacent to $z_u$.
	Note that~$z_u$ has exactly one neighbor~$z_u^k$ in~$Q^k_1$.
	When $z_u^k$ is not adjacent to $v$, any cycle through~$v$ suffices.
	Otherwise, we observe that $v$ has at least four neighbors in~$Q^k_1$, where one belongs to~$M$ and one may be $z_u^k$.
	For the~remaining two neighbors, we apply Theorem~\ref{theorem:alahmadi} to obtain a~matching~$M$ extending cycle $C_v = (v, \dots, z_v)$ where $z_v$ is is one of the~neighbors.
	By construction, $z_v \neq z_u^k$, yielding the~desired non-adjacent endpoints.

	Now suppose condition C2 holds, meaning~$M$ contains a~$u$-avoiding almost half-layer in direction $i$ where $v = u^i$, along with edges $u u^j$ and $v v^j$ in~$M$ for some~$j \neq i$.
	The~graph $Q_n - \set{u, u^j}$ contains two almost half-layers but no complete half-layer.
	Applying Theorem~\ref{theorem:skrekovski-2} to this graph with matching~$M \setminus \set{u u^j}$, we obtain a~cycle $C = (v, v^j, \dots, z_v)$ where $z_v = v^k$ for some $k$ distinct from both $i$ and $j$.
	Setting $z_u = u^j$, we can verify these endpoints are not adjacent: $z_u = (v^i)^j$ while $z_v = v^k$, and these vertices differ in multiple coordinates.
\end{proof}

The~proof of Theorem~\ref{theorem:hamilton-path} builds upon techniques used in proving Theorem~\ref{theorem:hamilton-cycle}, heavily utilizing the~properties of cube sequences established in Lemma~\ref{lemma:path-sequence}.

\begin{proof}[Proof of Theorem~\ref{theorem:hamilton-path}]
	For $n = d$, the~theorem follows directly from Conjecture~\ref{conjecture:main-path}, so we assume $n > d$ and, without loss of generality, that~$M$ spans only the~last $d$ directions.
	Throughout the~proof, we refer to cubes spanning the~last $d$ dimensions when discussing canonical copies of~$Q_d$ in~$Q_n$.
	We also use $k = |V(Q_{n - d})|$.

	The~proof strategy follows a~similar structure to Theorem~\ref{theorem:hamilton-cycle}.
	We locate appropriate paths in~$Q_{n-d}$, construct Hamilton paths within the~subcubes~$Q_d$, and connect these paths along the~structure in~$Q_{n-d}$.
	First, we introduce a~key construction that will be used throughout the~proof.

	The~proof proceeds by considering two main cases based on whether the~vertices lie in the~same canonical copy of~$Q_d$, with several subcases for each situation.

	\textbf{Case 1:} $u^L \neq v^L$ (vertices in different canonical copies)

	\textbf{Case 1.1:} $u^L$ and $v^L$ have different parity.

	Let $Y = (u^L,...,v^L)$ be a~Hamilton path in~$Q_{n-d}$ from $u^L$ to $v^L$.
	Since no C-conditions hold, the~induced sequence of matchings by $Y$ cannot contain a~sequence of half-layers covering both $u$ and $v$.
	Lemma~\ref{lemma:path-sequence} then guarantees the~existence of our desired Hamilton path.

	\textbf{Case 1.2:} $u^L$ and $v^L$ have same parity.

	A~direct Hamilton path between $u^L$ and $v^L$ in~$Q_{n-d}$ cannot exist due to their matching parity.
	We instead identify a~vertex $z$, distinct from $u^L$ and $v^L$ and of opposite parity to them, that we can avoid in our path construction.
	By known result from~\cite{lewinter97}, specifically Corollary~4, we get a~Hamilton path $Y_z$ between $u^L$ and $v^L$ on~$Q_{n - d} - \set{z}$.

	This vertex $z$ must be chosen such that the~induced sequence of cubes from the~resulting Hamilton path $Y_z$ on~$Q_{n-d} - \set{z}$ contains no sequence of half-layers covering both $u^R$ and $v^R$.
	Such a~choice must exist - otherwise, each canonical copy of~$Q_d$ would contain a~half-layer of the~same global parity, implying~$M$ contains a~half-layer covering both $u$ and $v$, contradicting our assumptions.

	Given this Hamilton path $Y_z$ and its induced sequence $G_z$ with matching~$M_z$, Lemma~\ref{lemma:path-sequence} provides a~Hamilton path $Y'_z$.
	To incorporate~$Q_n[z]$, we select any neighbor $z'$ of $z$ in~$Q_{n-d}$ and examine the~portion of $Y'_z$ within~$Q_n[z']$.
	From Lemma~\ref{lemma:path-sequence}, the~path construction provides at least~$2^{d-3}$ distinct choices for endpoint selection in the~cube~$Q_n[z']$, while maintaining the~existence of Hamilton path $Y'_z$ on the~induced sequence of cubes.
	Since $2^{d-3} \geq 2$ for $d \geq 5$, we can apply Lemma~\ref{lemma:cube-attachment} to complete this case by attaching the~subcube~$Q_n[z]$ to $Y'_z$ while extending~$M[z]$.

	\textbf{Case 2:} $u^L = v^L$ (vertices in same canonical copy), let~$M^\star = M[u^L]$ and~$Q^\star = Q_n[u^L]$ denote the~matching and cube containing both vertices.

	\textbf{Case 2.1:} $n > d + 1$

	\textbf{Case 2.1.1:} Hamilton path $Y^\star$ extending~$M^\star$ between $u^R$ and $v^R$ exists in~$Q^\star$.

	Consider a~cycle $C = (u^L, c_2, \dots, c_k)$ on~$Q_{n-d}$ and let $Y = (c_2, \dots, c_k)$.
	For any edge $wz \in Y^\star \setminus M^\star$, examine its corresponding vertices $w'$ in~$Q_n[c_2]$ and $z'$ in~$Q_n[c_k]$.
	If Lemma~\ref{lemma:path-sequence} provides a~Hamilton path $Y'$ extending the~induced matching~$M'$ between $w'$ and $z'$, we construct our desired path as $Y^\star \setminus \set{ w z } \cup \set{ w w', z z' } \cup Y'$.

	If no such edge exists, every pair $w'$ and $z'$ must be covered by a~sequence of half-layers in direction $i$ contained in~$M'$.
	Since $Y$ has odd length, $w'$ and $z'$ are covered by half-layers of same local parity, forcing edge $wz$ to share direction $i$.
	This implies all edges in~$Y^\star \setminus M^\star$ share direction~$i$.

	Consequently,~$M^\star$ must be perfect or almost perfect, as any uncovered inner vertex of $Y^\star$ would yield two incident edges in different directions.
	Moreover,~$M'$ must contain a~sequence of layers, not just half-layers, as $Y^\star \setminus M^\star$ cannot be covered by a~single half-layer.

	Thus~$M$ is either perfect or almost perfect with only endpoints uncovered.
	In the~almost perfect case, apply Theorem~\ref{theorem:perfect-matching} to~$M \cup \set{ u v }$ to obtain matching $N$, yielding Hamilton path~$M \cup N$.
	In the~perfect case, apply Theorem~\ref{theorem:skrekovski}, noting no C-conditions hold.

	\textbf{Case 2.1.2}: Hamilton path between $u$ and $v$ extending~$M^\star$ does not exist in~$Q^\star$.

	Since $uv \notin M$, condition C3 doesn't hold for~$M^\star$, and either condition C1 or C2 must hold for~$M^\star$.
	Applying Lemma~\ref{lemma:disjoint-paths} to $u$, $v$ and~$M^\star$ yields two vertex-disjoint paths~$Z_u$ and~$Z_v$ terminating at vertices $z_u$ and $z_v$ respectively, together covering all vertices in~$Q^\star$ and extending~$M^\star$.

	Following the~approach from the~previous case, we consider a~cycle $C = (u^L, c_2, \dots, c_k)$ on~$Q_{n-d}$ and examine the~induced sequence of cubes by $Y = (c_2, \dots, c_k)$.
	If there exists a~Hamilton path~$Y'$ extending the~induced matching~$M'$ between corresponding $z'_u \in Q_n[c_2]$ and $z'_v \in Q_n[c_k]$, then we are done since $Z_u \cup Z_v \cup Y' \cup \set{ z_u z'_u, z_v z'_v }$ is the~needed Hamilton path.

	Suppose that there is no such path.
	Then the~sequence of cubes induced by $Y$ contains a~sequence of half-layers covering $z'_u$ and $z'_v$ in a~direction $i$.
	Since for~$M^\star$ holds one of the~conditions~C1 and~C2 we know that~$M^\star$ contains a~half-layer or an~almost half-layer in a~direction~$j$ and both~$z'_u$ and~$z'_v$ are not covered by it.
	We say that $i \neq j$, since otherwise the~sequence of half-layers contained in~$M'$ has the~same direction and global parity as the~half-layer or an~almost half-layer in~$M^\star$, violating the~C-conditions globally.

	So $i \neq j$, and we can uniquely extend this sequence to layers and similarly extend~$M^\star$ to a~perfect matching both containing obstructions in different directions.
	Therefore we can apply Theorem~\ref{theorem:skrekovski} to obtain our desired path.

	\textbf{Case 2.2:} $n = d + 1$, let~$Q^\dagger$ denote the~canonical copy of~$Q_d$ adjacent to~$Q^\star$ in~$Q_n$, and let~$M^\dagger$ be the~restriction of~$M$ to~$Q^\dagger$.

	\textbf{Case 2.2.1:} Hamilton path $Y^\star$ extending~$M^\star$ between $u^R$ and $v^R$ exists in~$Q^\star$.

	From Conjecture~\ref{conjecture:main-path-two} we get existence of two distinct Hamilton paths extending~$M^\star$ between $u^R$ and $v^R$ in~$Q^\star$ and we are free to use Lemma~\ref{lemma:cube-attachment} to obtain the~final result.

	\textbf{Case 2.2.2:} Hamilton path between $u^R$ and $v^R$ does not exist in~$Q^\star$.
	We follow the~approach from Case 2.1.2.
	Lemma~\ref{lemma:disjoint-paths} provides two vertex-disjoint paths~$Z_u$ starting from~$u$ and~$Z_v$ starting from~$v$ in~$Q^\star$, terminating at non-adjacent vertices~$z_u$ and~$z_v$, with corresponding vertices~$z_u'$ and~$z_v'$ in~$Q^\dagger$.
	Since~$z_u$ and~$z_v$ are not neighbors, $z'_u z'_v \notin M^\dagger$.

	If~$M^\dagger$ contains no half-layer, Conjecture~\ref{conjecture:main-path} guarantees a~Hamilton path between~$z_u'$ and~$z_v'$ extending~$M^\dagger$, since they are not neighbors.
	Otherwise, both~$M^\star$ and~$M^\dagger$ can be uniquely extended to perfect matchings and no C-conditions hold globally for~$M$, by the~discussion from Case 2.1.2.
	Allowing us to apply Theorem~\ref{theorem:skrekovski} to obtain the~desired path.
\end{proof}

\begin{proof}[Proof of Theorem~\ref{theorem:hamilton-path-5}]
	The~theorem directly follows from Observation \ref{observation:main-path-verification} and Theorem \ref{theorem:hamilton-path}.
\end{proof}

\bibliographystyle{splncs04}
\bibliography{cube}

\begin{thebibliography}{10}
\providecommand{\url}[1]{\texttt{#1}}
\providecommand{\urlprefix}{URL }
\providecommand{\doi}[1]{https://doi.org/#1}

\bibitem{thomassen2015extending}
Alahmadi, A., Aldred, R., Alkenani, A., Hijazi, R., Sol{\'e}, P., Thomassen,
  C.: {Extending a perfect matching to a Hamiltonian cycle}. Discrete
  Mathematics \& Theoretical Computer Science  \textbf{17} (2015)

\bibitem{alahmadi2015extending}
Alahmadi, A., Aldred, R.E., Alkenani, A., Hijazi, R., Sol{\'e}, P., Thomassen,
  C.: {Extending a perfect matching to a Hamiltonian cycle}. Discrete
  Mathematics and Theoretical Computer Science (Online Edition)
  \textbf{17}(1),  241--254 (2015)

\bibitem{biere2024cadical}
Biere, A., Faller, T., Fazekas, K., Fleury, M., Froleyks, N., Pollitt, F.:
  {CaDiCaL 2.0}. In: Gurfinkel, A., Ganesh, V. (eds.) Computer Aided
  Verification - 36th International Conference, {CAV } 2024, Montreal, QC,
  Canada, July 24-27, 2024, Proceedings, Part {I}. Lecture Notes in Computer
  Science, vol. 14681, pp. 133--152. Springer (2024).
  \doi{10.1007/978-3-031-65627-9_7}

\bibitem{bright2016mathcheck}
Bright, C., Ganesh, V., Heinle, A., Kotsireas, I., Nejati, S., Czarnecki, K.:
  \textsc{MathCheck2}: A {SAT}+{CAS} Verifier for Combinatorial Conjectures, p.
  117–133. Springer International Publishing (2016).
  \doi{10.1007/978-3-319-45641-6_9}

\bibitem{skrekovski2009gray}
Dimitrov, D., Dvo{\v{r}}{\'a}k, T., Gregor, P., {\v{S}}~krekovski, R.: Gray
  codes avoiding matchings. Discrete Mathematics \& Theoretical Computer
  Science  \textbf{11},  123--147 (2009)

\bibitem{dvorak2016}
Dvo{\v{r}}{\'a}k, T., Fink, J.: Gray codes extending quadratic matchings. J
  Graph Theory  (2018). \doi{10.1002/jgt.22371}

\bibitem{dvovrak2007hamiltonian}
Dvo{\v{r}}{\'a}k, T., Gregor, P.: Hamiltonian paths with prescribed edges in
  hypercubes. Discrete Mathematics  \textbf{307}(16),  1982--1998 (2007)

\bibitem{dvorak2008faulty}
Dvo{\v{r}}{\'a}k, T., Gregor, P.: Partitions of faulty hypercubes into paths
  with prescribed endvertices. SIAM Journal on Discrete Mathematics
  \textbf{22}(4),  1448--1461 (2008)

\bibitem{Dvorak_cycle_prescribed}
Dvo\v{r}\'{a}k, T.: Hamiltonian cycles with prescribed edges in hypercubes.
  SIAM J. Discret. Math.  \textbf{19}(1),  135--144 (2005)

\bibitem{kreweras1}
Fink, J.: Perfect matchings extend to {Hamilton} cycles in hypercubes. J. Comb.
  Theory, Ser. B  \textbf{97}(6),  1074--1076 (2007)

\bibitem{2factor}
Fink, J.: Matchings extend into 2-factors in hypercubes. Combinatorica
  \textbf{39},  77--84 (2019). \doi{10.1007/s00493-017-3731-8}

\bibitem{fink2024matchings}
Fink, J., M{\"u}tze, T.: Matchings in hypercubes extend to long cycles. In:
  International Workshop on Combinatorial Algorithms. pp. 14--27. Springer
  (2024)

\bibitem{gregor2009perfect}
Gregor, P.: {Perfect matchings extending on subcubes to Hamiltonian cycles of
  hypercubes}. Discrete Mathematics  \textbf{309}(6),  1711--1713 (2009)

\bibitem{gregor2017}
Gregor, P., Novotn{\'y}, T., {\v{S}}krekovski, R.: Extending perfect matchings
  to gray codes with prescribed ends. The Electronic Journal of Combinatorics
  \textbf{25}(2),  2--56 (2018)

\bibitem{Gros}
Gros, L.: Th\'{e}orie du Baguenodier. Aim\'{e} Vingtrinier, Lyon (1872)

\bibitem{knuth2005art}
Knuth, D.E.: The Art of Computer Programming, Volume 4, Fascicle 2: Generating
  All Tuples and Permutations. Addison-Wesley Professional (2005)

\bibitem{Kreweras}
Kreweras, G.: {Matchings and Hamiltonian cycles on hypercubes}. Bull. Inst.
  Combin. Appl.  \textbf{16},  87--91 (1996)

\bibitem{lewinter97}
Lewinter, M., Widulski, W.: Hyper-hamilton laceable and caterpillar-spannable
  product graphs. Computers \& Mathematics with Applications  \textbf{34}(11),
  99–104 (1997). \doi{10.1016/s0898-1221(97)00223-x}

\bibitem{Ruskey}
Ruskey, F., Savage, C.: {Hamilton Cycles that Extend Transposition Matchings in
  Cayley Graphs of $S_n$}. SIAM Journal on Discrete Mathematics  \textbf{6}(1),
   152--166 (1993)

\bibitem{Savage}
Savage, C.: {A survey of combinatorial Gray codes}. SIAM Review
  \textbf{39}(4),  605--629 (1997)

\bibitem{simmons1977almost}
Simmons, G.J.: Almost all n-dimensional rectangular lattices are
  hamilton-laceable. In: Presented at the 9th Southeastern Conf. on
  Combinatoric Graph Theory and Computing (1977)

\bibitem{vandenbussche2013extensions}
Vandenbussche, J., West, D.B.: Extensions to $2 $-factors in bipartite graphs.
  The Electronic Journal of Combinatorics  \textbf{20}(3),  1--10 (2013)

\bibitem{wang24}
Wang, S., Wang, F.: A kind of matchings extend to hamiltonian cycles in
  hypercubes. RAIRO-Operations Research  \textbf{58}(6),  5237--5254 (2024)

\bibitem{xu2013topological}
Xu, J.: Topological structure and analysis of interconnection networks, vol.~7.
  Springer Science \& Business Media (2013)

\end{thebibliography}

\newpage
\appendix

\section{Relations between statements}
\label{sec:relation}

In this section, we prove that Conjecture \ref{conjecture:main-path} implies Conjecture \ref{conjecture:maximal} and Conjecture~\ref{conjecture:main-cycle} and the~last one implies Ruskey-Savage Conjecture~\ref{conjecture:ruskey-savage}.
Then, we show that Conjecture \ref{conjecture:main-path} cannot be generalized for matchings of $B(Q_n)$ using the~construction of \cite{dvorak2016}.

\begin{observation}
	\label{observation:maximal-implies}
	Conjecture~\ref{conjecture:main-path} implies Conjecture~\ref{conjecture:maximal}.
\end{observation}
\begin{proof}
	Let~$M$ be a~maximal matching.
	We prove the~equivalence between a~C-condition holding and $M' = (M \cup \set{u^M v^M}) \setminus \set{u u^M, v v^M}$ containing a~half-layer.

	First, suppose one of the~C-conditions holds.
	If C1 holds, then~$M$ contains a~half-layer covering~$u$ and $v$.
	Since~$M$ is maximal, it contains both half-layers in this direction, and the~operation creating~$M'$ preserves one of them.
	When C2 holds, $M$ contains a~$u$-avoiding almost half-layer in direction $i$ where $v = u^i$, and adding the~edge~$u^M v^M$ completes the~another almost half-layer in~$M'$.

	For the~converse, suppose~$M'$ contains a~half-layer.
	If~$M$ already contains a~half-layer, then by maximality of~$M$, it contains both half-layers in this direction.
	Since~$M'$ still contains a~half-layer after removing edges~$u u^M$ and $v v^M$, vertices~$u$ and $v$ must be covered by the~same half-layer, satisfying C1.
	Otherwise, the~half-layer in~$M'$ must have been created by adding~$u^M v^M$, implying~$M$ contains an~almost half-layer satisfying C2.
	Since~$u v \notin M$ by assumption, C3 cannot occur.
\end{proof}

\begin{observation}
	\label{observation:conjecture-path-implies-cycle}
	Conjecture~\ref{conjecture:main-path} implies Conjecture~\ref{conjecture:main-cycle}, for dimension $d \geq 5$.
\end{observation}
\begin{proof}
	Since one of the~vertices $u$ and $v$ is not covered, none of the~C-conditions may hold and thus there exists a~Hamilton path extending matching~$M$ between $u$ and $v$.
\end{proof}

\begin{observation}
	\label{observation:main-cycle-implies-ruskey-savage}
	Conjecture~\ref{conjecture:main-cycle} implies Conjecture~\ref{conjecture:ruskey-savage}.
\end{observation}
\begin{proof}
	For a~perfect matching, the~observation follows directly from Theorem~\ref{theorem:perfect-matching}.
	If the~matching is not perfect, there exists a~vertex $u$ avoided by~$M$.
	Let $v$ be a~neighbor of $u$.
	We immediately see that~$M$, $u$ and $v$ fulfill Conjecture~\ref{conjecture:main-cycle} criteria.
	Thus, a~Hamilton path $Z$ exists between~$u$ and~$v$.
	This path~$Z$, along with the~edge $u v$, forms a~Hamilton cycle extending~$M$, and Conjecture~\ref{conjecture:ruskey-savage} follows.
\end{proof}

%\section{Matching not extendable to Hamilton path in $B(Q_n)$}
%\label{sec:matching-not-extendable}

For a~vertex $v \in V(Q_n)$ we denote $\chi(v)$ to be $1$ for even vertices and $-1$ for odd vertices.
We also need the~following proposition:

\begin{proposition}[Dvo\v{r}ák and Gregor~\cite{dvorak2008faulty}, Proposition 3.1]
	\label{proposition:balanced}
	Let $H$ be a~subgraph of~$Q_n$.
	Let $\set{Z_{u_i v_i}}^m_{i=1}$ be paths in $H$ that are pairwise vertex-disjoint, between endpoints $u_i$ and $v_i$, such that they together cover all the~vertices of $H$, then

	\begin{equation}
		\label{eq:balanced}
		\frac{1}{2} \sum^m_{i=1} (\chi(u_i) + \chi(v_i)) = \sum_{u \in V(H)} \chi(u).
	\end{equation}
\end{proposition}

Note that in the~proposition statement $u_i$ and $v_i$ may be the~same vertex, which means we can have a~paths consisting of single vertex.

\begin{observation}
	\label{observation:BQn-not-extendable}
	For every $n \ge 9$, there exists a~matching~$M$ of $B(Q_n)$ that cannot be extended to a~Hamilton path using only edges from~$Q_n$.
\end{observation}
\begin{proof}
	This proof follows the~approach of Theorem~1.2 in~\cite{dvorak2016}.
	We begin with inequality that holds for all $m \geq 8$:
	\begin{equation}
		\label{eq:binom-inequality}
		\binom{m}{\floor{\frac{m}{2}}} + 3
		\leq \sum^{\floor{ \frac{m}{2} } - 1}_{i=0} \binom{m}{i}
		\leq \sum^{m}_{i = \floor{ m/2 } + 1} \binom{m}{i}.
	\end{equation}

	Let~$G_1$ and~$G_2$ be subgraphs of~$B(Q_n)$ induced by vertices with first coordinate~$0$ and~$1$ respectively.
	Note that for any vertex~$a \in V(G_1)$, the~vertex~$a^1 \in V(G_2)$.
	Let~$h := \floor{ \frac{n-1}{2} }$ and~$k~:=~\binom{n - 1}{h}$.
	We partition~$G_1$ and~$G_2$ into three layers:
	\begin{align*}
		\sslow{}^0 & := \set{ v \in V(G_1) \mid w(v) < h }, & \sslow{}^1 & := \set{ v^1 \mid v \in \sslow{}^0 }, \\
		\ssmid{}^0 & := \set{ v \in V(G_1) \mid w(v) = h }, & \ssmid{}^1 & := \set{ v^1 \mid v \in \ssmid{}^0 }, \\
		\ssup{}^0  & := \set{ v \in V(G_1) \mid w(v) > h }, & \ssup{}^1  & := \set{ v^1 \mid v \in \ssup{}^0 }.
	\end{align*}

	Let $\set{x_1, \dots, x_k}$ enumerate the~vertices in $\ssmid{}^0$ and put $p := \chi(x_1)$.
	From~(\ref{eq:binom-inequality}) with $m = n - 1$, we have $|\ssup{}^0| \geq |\sslow{}^0| \geq |\ssmid{}^0| + 3$.
	This allows us to select distinct vertices $y_1, \dots, y_{k+3} \in \sslow{}^0 \cup \sslow{}^1$ and $z_1, \dots, z_{k + 3} \in \ssup{}^0 \cup \ssup{}^1$ satisfying $\chi(z_i) = p = -\chi(y_i)$ for all $i \in [k + 3]$.
	We construct matching~$M$ in~$B(Q_n)$ as $\set{x_i y_i}^k_{i=1} \cup \set{x^1_i z_i}^k_{i=1} \cup \set{y_{k + i} z_{k + i}}^3_{i=1}$.

	Suppose for contradiction that there exists a~Hamilton path $Z$, that must be between vertices of opposite parity, that extends~$M$ using only edges from~$Q_n$.
	Let $\sslm{}_n$ be the~subgraph of~$Q_n$ induced by vertices $\sslow{}^0 \cup \sslow{}^1 \cup \ssmid{}^0 \cup \ssmid{}^1$.
	The~edges $(E(Z) \setminus M) \cap E(\sslm{}_n)$ form a~collection~$\mathcal{P}$ of vertex-disjoint paths covering~$\sslm{}_n$, with endpoints in $(V(M) \cup \set{u, v}) \cap V(\sslm{}_n)$.
	Also $\sslow{}^0$ contains the~same number of even vertices as $\sslow{}^1$ of odd vertices and vice-versa.
	That also applies to $\ssmid{}^0$ and $\ssmid{}^1$ and Proposition~\ref{proposition:balanced} yields
	\begin{equation}
		\label{eq:balanced-contradiction}
		\sum_{u \in V(\sslm{}_n)} \chi(u) = 0 \underset{\text{(\ref{eq:balanced})}}{=} \sum_{Z_{u v} \in \mathcal{P}} \overbrace{\chi(u) + \chi(v)}^{\text{($\star$)}}.
	\end{equation}

	Each path endpoint contributes either $\chi(u)$ or $2\chi(u)$ to ($\star$), the~latter occurring for single-vertex paths.
	By construction of~$M$, we have at least $2k + 3$ vertices of parity $-p$ and at most~$k + 1$ vertices of parity $p$ among the~endpoints, implying
	\[
		\left|\sum_{Z_{u v} \in \mathcal{P}} \chi(u) + \chi(v)\right| > 0,
	\]
	contradicting (\ref{eq:balanced-contradiction}).
\end{proof}

\section{Computer verification}
\label{sec:computer}

We establish the~validity of Conjecture~\ref{conjecture:main-cycle} for dimensions $2 \leq n \leq 5$, along with Conjecture~\ref{conjecture:main-path-two} for $n = 5$ through computer-assisted verification techniques inspired by the~\smathcheck{} methodology~\cite{bright2016mathcheck}.
The~source code for this section is available at \url{https://gitlab.com/jirka.fink/matching_of_five_directions}.

Our verification strategy centers on systematically enumerating feasible matchings of~$Q_n$ and verifying the~existence of Hamilton paths between vertices of differing parity.
To accomplish this, we develop a~\ssat{} representation of the~matching constraints and employ the~\scadical{} solver~\cite{biere2024cadical} for matching generation.
After each matching is generated and verified, we exclude both it and its isomorphic variants from subsequent iterations, ensuring comprehensive coverage of the~solution space.

Given the~vast number of possible matchings, we optimize our analysis by focusing primarily on maximal matchings.
This approach is justified since any non-maximal matching extends to a~maximal configuration and inherits its path extension properties, except of the~matchings subject to the~C-conditions, that will be handled separately for Conjecture~\ref{conjecture:main-path-two}.
For Conjecture~\ref{conjecture:main-cycle}, we further optimize by forcing vertex $\bm{0}$ to remain uncovered during maximal matching generation, as the~conjecture requires at least one avoided vertex.

The~\ssat{} formulation for matching identification in dimension $n$ introduces variables $x_e$ for each edge $e \in E(Q_n)$, indicating membership in matching~$M$.
To enforce matching properties, we specify that no two incident edges can be selected:

\[
	\bigwedge_{\substack{e,f \in E(Q_n) \\ e \ \text{incident to}\ f}} \neg x_e \lor \neg x_f.
\]

We introduce vertex coverage variables $x_v$ for each $v \in V(Q_n)$, linking them to incident edge selection:
\[
	\bigwedge_{v \in V(Q_n)} \left( \bigvee_{\substack{e \in E(Q_n) \\ e\ \text{incident to}\ v}} x_e \right) \lor \neg x_v.
\]

We also add there the~opposite implication, that is, presence of an~edge $u v$ implies that the~vertices $u$ and $v$ are covered:
\[
	\bigwedge_{u v \in E(Q_n)} \left( \neg x_{uv} \lor x_u \right) \wedge \left( \neg x_{uv} \lor x_v \right).
\]

To ensure the~maximality condition at least one vertex in each adjacent pair must be covered:
\[
	\bigwedge_{\substack{u, v \in V(Q_n) \\ u\ \text{adjacent to}\ v}} x_u \lor x_v.
\]

To exclude matchings already verified, we introduce a~set $\mathcal{F}$, representing all previously verified matchings.
\[
	\forall F \in \mathcal{F}: \bigvee_{e \in E(Q_n) \setminus F} x_e.
\]

This means that the~newly generated matching has to contain at least one edge that is not present in $F, \forall F \in \mathcal{F}$.
Here as previously specified we also include in the~set $\mathcal{F}$ all isomorphic matchings to the~newly generated matching $F$.

For Conjecture~\ref{conjecture:main-cycle}, we explicitly avoid the~vertex $\bm{0}$:
\[
	\neg x_{\bm{0}}.
\]

The~verification process systematically examines terminal vertex pairs for each maximal matching.
For Conjecture~\ref{conjecture:main-cycle}, we pair vertex $\bm{0}$ with each odd vertex, while for Conjecture~\ref{conjecture:main-path-two}, we consider all opposite-parity vertex pairs not subject to C-conditions.
This direct verification establishes the~following observation:

\begin{observation}
	\label{observation:main-cycle-verification}
	Conjecture~\ref{conjecture:main-cycle} holds for $2 \leq d \leq 5$.
\end{observation}

For Conjecture~\ref{conjecture:main-path-two}, we must also verify matchings~$M$ and terminal pairs $u, v$ where adding a~single edge $e$ would cause~$M \cup \set{e}$ to violate C-conditions.
We identify such matchings through a~two-step process: First, we start with maximal matchings that violate C-conditions and iteratively remove edges until reaching minimal violating configurations.
Then, from each minimal configuration, we remove one edge at a~time to obtain matchings that are exactly one edge away from violating C-conditions.
For these matchings, we verify the~existence of two distinct Hamilton paths between each valid terminal pair.
This approach completes our proof for $n = 5$, as any matching not violating C-conditions either extends to a~C-condition-free maximal matching or extends to a~matching that is one edge away from violating C-conditions.
These verifications establish:

\begin{observation}
	\label{observation:main-path-verification}
	Conjecture~\ref{conjecture:main-path-two} holds for $d = 5$.
\end{observation}

Hamilton paths are found using the~backtracking algorithm from~\cite{fink2024matchings}.
For finding a~second path, we iteratively force the~usage of edges absent from the~first path during backtracking initialization until a~second valid path is discovered.

Conjecture~\ref{conjecture:main-path-two} fails for $n=4$, as we encounter the~obstructions previously identified by Gregor, Novotn\'{y} and \v{S}krekovski~\cite{gregor2017}.
Additionally, we discovered a~novel obstruction specific to maximal (but not perfect) matchings, illustrated in Figure~\ref{figure:forbidden-matchings}.

Initial verification attempts on an~Apple MacBook M3 Max proved intractable since the~generation of maximal matchings did not finish even after 3 days.
We achieved significant speedup by partitioning the~search space: since maximal matchings either cover all vertices or avoid specific non-adjacent vertex combinations, we generated these combinations up to isomorphism.
There were 159 such combinations, and we verified them in parallel.

The~optimized algorithm reduced maximal matching generation for $n = 5$ to 6 hours, producing 16,459 non-isomorphic and 59,457,409 isomorphic matchings, with Hamilton path verification completing in 2 minutes.
A~specialized handcrafted algorithm, inspired by~\cite{fink2024matchings}, further reduced matching generation time to 3 minutes.
The~subsequent verification of the~cases one edge away from violating C-conditions takes 20 minutes.
For $n = 4$, complete verification executes in under 10 seconds.

\end{document}